\documentclass[10pt,reqno]{amsart}
\usepackage{amsmath,amsthm,
amsopn,amssymb} 

\DeclareMathOperator{\ham}{ham} \DeclareMathOperator{\ann}{ann}
 
 \DeclareMathOperator{\Pann}{Pann}
\DeclareMathOperator{\ad}{ad} 

 \DeclareMathOperator{\Endo}{End}
 \DeclareMathOperator{\PZ}{PZ}
\DeclareMathOperator{\sing}{sing}
\newcommand{\ov}{\overline}

\newcommand{\Z}{{\mathbb Z}}
\newcommand{\C}{{\mathbb C}}

\newcommand{\F}{{\mathbb F}}
\newcommand{\Q}{{\mathbb Q}}
\newcommand{\g}{{\mathfrak g}}
\newcommand{\h}{{\mathfrak h}}
\newcommand{\s}{{\mathfrak s}}
\newcommand{\m}{{\mathfrak m}}

\newcommand{\leftexp}[2]{{\vphantom{#2}}^{#1\hspace{-.5mm}}{#2}}

\theoremstyle{plain}
\newtheorem{theorem}{Theorem}[section]
\newtheorem{lemma}[theorem]{Lemma}

\newtheorem{prop}[theorem]{Proposition}

\theoremstyle{definition}
\newtheorem{defn}[theorem]{Definition}
\newtheorem{defns}[theorem]{Definitions}
\newtheorem{rmk}[theorem]{Remark}
\newtheorem{rmks}[theorem]{Remarks}
\newtheorem{notn}[theorem]{Notation}
\newtheorem{example}[theorem]{Example}

\newcommand{\J}{\ov{J}}

\numberwithin{equation}{section}

\begin{document}
\bibliographystyle{amsalpha}
\author[D. A. Jordan]{David A. Jordan}

\address
{Department of Pure Mathematics\\
University of  Sheffield\\
Hicks Building\\
Sheffield S3~7RH\\
UK}

\email{d.a.jordan@sheffield.ac.uk}

\date{}

\title
{Finite-dimensional simple Poisson modules}


\begin{abstract}We prove a result that can be applied to determine the finite-dimensional
simple Poisson modules over a Poisson algebra and apply it to
numerous examples. In the discussion of the examples, the emphasis
is on the correspondence with the finite-dimensional simple modules
over deformations and on the behaviour of finite-dimensional simple
Poisson modules on the passage from a Poisson algebra to the Poisson
subalgebra of invariants for the action of a finite group of Poisson
automorphisms.
\end{abstract}

 \maketitle

\section{Introduction}
There are several known situations where aspects of a Poisson
algebra $A$ and of a noncommutative deformation $T$ of $A$ reflect
and enlighten each other. Examples include Poisson primitive ideals
and primitive ideals
\cite{hodlev,joseph,oh1,ChoOh,goodlaun,goodearlDMP,goodearlDMq} and
Poisson (co)homology and Hochschild (co)homology
\cite{Kasselpaper,alevlambre,alevfoissy}. There is rarely a formal
mechanism for passing between the two, an exception being the
Dixmier map, as specified in \cite{kab}, which, following
\cite{mathieu}, is known to give a homeomorphism between the prime
spectrum of the enveloping algebra $U(\g)$ of a solvable Lie algebra
and the Poisson prime spectrum of the symmetric algebra $S(\g)$,
under the Kirillov-Kostant bracket \cite[III.5.5]{BGl}. The prime
ideals of $U(\g)$ are completely prime but if $T$ has
finite-dimensional simple modules of dimension $>1$ then their
annihilators are primitive ideals that are not completely prime and
have no clear counterpart among the prime Poisson ideals of $A$.
Although such ideals are not reflected in the Poisson structure, we
shall see that the simple modules that they annihilate often are. We
shall present a result giving a method for the determination of the
finite-dimensional simple Poisson modules over any affine Poisson
algebra. Here the term
 \textit{Poisson module} is used in the sense of
Farkas \cite{farkas} and Oh \cite{oh}. If $J$ is a Poisson maximal
ideal of $A$ then $J/J^2$ has a natural Lie algebra structure and we
show that there is a  bijection, preserving dimension, between the
isomorphism classes of finite-dimensional simple Poisson $A$-modules
and pairs $(J,\widehat{M})$ where $J$ is a Poisson maximal ideal of
$A$ and $\widehat{M}$ is an isomorphism class of finite-dimensional
modules over the Lie algebra $J/J^2$. In this bijection, the simple
Poisson modules in a class corresponding to the pair
$(J,\widehat{M})$ are annihilated by $J$. The proof does not require
any deep ideas and is essentially a matter of lifting various module
structures, associative and Lie, and checking compatibility. It is
constructive inasmuch that if the Poisson maximal ideals $J$ of $A$
can be identified and the representation theory of the Lie algebras
$J/J^2$ is known then the finite-dimensional simple Poisson
$A$-modules can be computed. As we shall see, this works very neatly
in many examples.

We illustrate the method by computing the finite-dimensional simple
Poisson modules for numerous examples and for many of these we
compare the results with the classification of finite-dimensional
simple modules over a noncommutative deformation of $A$, in an
algebraic sense to be specified in \ref{quantdeform}, for which this
classification is known. Let us say that a $\C$-algebra $T$ is
\textit{$t$-homogeneous} if, for, each positive integer $d$, there
are, up to isomorphism, precisely $t$ simple (left) $T$-modules of
dimension $d$ and that a Poisson algebra is \textit{$t$-homogeneous}
if it has the analogous Poisson property. The enveloping algebra
$U(sl_2)\simeq U(so_3)$ is well-known to be $1$-homogeneous, whereas
the quantized enveloping algebra $U_q(sl_2)$ is $2$-homogeneous (and
has a $4$-homogeneous variant) and the cyclically $q$-deformed
algebra $U_q^\prime (so_3)$ of \cite{mp1,mp2} is $5$-homogeneous.
Each of $U(sl_2)$, $U_q(sl_2)$ and $U_q^\prime (so_3)$ is a
deformation of a Poisson algebra $A$ and we shall see that these
Poisson algebras are $t$-homogeneous for $t=1, 2, 5$, respectively.
Here $t$ is, in each case, the number of Poisson maximal ideals of
$A$, each of which corresponds to a singularity of a surface, so, in
particular, this gives a geometric explanation for the mysterious
occurrence of the number five in the classification of simple
$U_q^\prime (so_3)$-modules in \cite{mp1,mp2}. The Poisson
$t$-homogeneity of these examples is due to the fact that
$J/J^2\simeq sl_2$ for each maximal Poisson ideal $J$ of $A$. We
include one example, \ref{Plocenv}, where the classification of the
finite-dimensional simple modules over the deformation was not
previously known but can be predicted using knowledge of the
finite-dimensional simple Poisson modules. Often, as in the examples
mentioned above, there is a bijection, preserving dimension, between
the isomorphism classes of finite-dimensional simple Poisson modules
and the isomorphism classes of finite-dimensional simple modules
over the corresponding deformation but, in other examples, we need
to replace a single deformation by a family of deformations for such
a bijection to exist.

Many of the Poisson algebras that we consider arise as the Poisson
algebra $B^G$ of invariants for the action of a finite group $G$ of
Poisson automorphisms of a Poisson algebra $B$. For example, the
algebra $U^\prime_q(so_3)$ has a factor that is the algebra of
invariants for an automorphism, of order $2$, of the coordinate ring
of the quantum torus and the analogous factor of the corresponding
Poisson algebra is the Poisson algebra of invariants for a Poisson
automorphism, of order $2$, of the coordinate ring of the torus. A
secondary theme of the paper is to explore aspects of the behaviour
of finite-dimensional simple Poisson modules on the passage from $B$
to $B^G$. Alev and Farkas \cite{alevfarkas} have pointed out that
$B^G$ need not be simple, as a Poisson algebra, when $B$ is simple.
In several such examples, including the one used in
\cite{alevfarkas}, there is one or more Poisson maximal ideal in
$B^G$, giving rise to finite-dimensional simple Poisson modules. We
shall observe various contrasts with the known theory for
associative algebras of invariants. In Example~\ref{KleinianA3up},
there is a $1$-homogeneous Poisson algebra $B$ for which every
finite simple Poisson $B$-module splits, over $B^G$, into a sum of
one-dimensional Poisson modules so that there are Poisson simple
$B$-modules of unbounded length over $B^G$. In contrast, in
Subsection~\ref{fourfour}, there is an example where both $B$ and
$B^G$ are $4$-homogeneous but, for each $d$, the four
$d$-dimensional simple Poisson $B$-modules each restrict to a simple
Poisson $B^G$-module, but only yielding two non-isomorphic simple
Poisson $B^G$-modules, and there are two $d$-dimensional simple
Poisson $B^G$-modules that are not summands of restrictions of
simple Poisson $B$-modules. In Example~\ref{Atwo}, $B$ is simple but
$B^G$ has a finite-dimensional Poisson module that is not
semisimple.

The paper is organised as follows. The basic ideas in the theory of
Poisson algebras and Poisson modules are set out in Section~2 and
the result on the determination of finite-dimensional simple Poisson
modules is in Section~3. Section~4 is subdivided and is concerned
with examples, most of which are either $\C[x,y,z]$, with a certain
type of Poisson bracket, or localizations or factors thereof. This
includes the Poisson algebras corresponding to $U^\prime_q(so_3)$
and $U_q(sl_2)$.  The final section deals with three interesting
examples, computed in \cite{alevfoissy}, of Poisson algebras of
invariants for the diagonal action of a Weyl group $W$ of a simple
Lie algebra $\g$ on the symmetric algebra
$S(\mathfrak{h}\oplus\mathfrak{h}^*)$ for a Cartan subalgebra $\h$
of $\g$. We show that, in each of the types $A_2, B_2$ and $G_2$,
the Poisson algebra of invariants is $1$-homogeneous.

The finite-dimensional simple Poisson modules associated to
$U^\prime_q(so_3)$ and $U_q(sl_2)$ were classified, using direct
methods, in the PhD thesis of N. Sasom \cite{nongthesis}.

We thank J. Alev for several helpful comments on Poisson algebras.

\section{Poisson algebras and modules}
\begin{defns}
Our base field will always be $\C$. By {\it Poisson algebra},
we mean a commutative $\C$-algebra $A$ with a bilinear product
$\{-,-\}:A\times A\rightarrow A$ such that $A$ is a Lie algebra
under $\{-,-\}$ and, for all $a\in A$, $\ham(a):\{a,-\}$ is a
$\C$-derivation of $A$, called a {\it Hamiltonian} derivation of $A$
(or a Hamiltonian vector field). Such a product $\{-,-\}$ will be
referred to as a {\it Poisson bracket} on $A$. By {\it affine
Poisson algebra}, we mean a Poisson algebra that is finitely
generated as a $\C$-algebra.

A subalgebra $B$ of a Poisson algebra $A$ is a {\it Poisson
subalgebra} of $A$ if $\{b,c\}\in B$ for all $b,c\in B$ and an ideal
$I$ of $A$ is a {\it Poisson ideal} if $\{i,a\}\in I$ for all $i\in
I$ and all $a\in A$, in which case $A/I$ is a Poisson algebra with
$\{a+I,b+I\}=\{a,b\}+I$ for all $a, b\in A$. A Poisson algebra is
\textit{simple} if it has no proper non-zero Poisson ideals. An
algebra homomorphism between two Poisson algebras is a
\textit{Poisson homomorphism} if it is also a Lie homomorphism and
is a \textit{Poisson isomorphism} if it is a bijective Poisson
homomorphism, in which case its inverse is also Poisson. Two Poisson
brackets on the same algebra $A$ are \textit{equivalent} if the
corresponding Poisson algebras are isomorphic. The  \textit{Poisson
centre} $\PZ(A)$ consists of those elements $a\in A$ such that
$\{a,b\}=0$  for all $b\in A$.
\end{defns}

\begin{rmks}\label{Pideals}
(i) Let $I$ and $J$ be Poisson ideals of a Poisson algebra $A$. As
$\ham(a)$ is a derivation of $A$ for all $a\in A$,  $IJ$ is a
Poisson ideal. Also $I$ and $J$ are Lie subalgebras of $A$ under
$\{-,-\}$ and, if $I\subseteq J$, $I$ is a Lie ideal of $J$ and
$J/I$ is a Lie algebra. In particular $J/J^2$ is always a Lie
algebra and will be denoted $\g(J)$.

(ii) There is a need to distinguish between the concepts of a
maximal ideal that is also a Poisson ideal, for which we use the
term {\it Poisson maximal ideal}, and a Poisson ideal that is
maximal in the lattice of Poisson ideals. An example of the latter,
but not the former, is $0$ in the simple Poisson algebra $\C[x,y]$,
with $\{x,y\}=1$.

(iii) If $A$ is affine then there is no need to distinguish between
the concepts of a prime ideal that is also a Poisson ideal and a
Poisson ideal that is prime in the lattice of Poisson ideals, in the
sense that for all Poisson ideals $I$ and $J$ of $A$, $IJ\subseteq
P$ implies $I\subseteq P$ or $J\subseteq P$. This is a consequence
of the well-known fact \cite[Lemma 3.3.3]{dix} that if $\delta$ is a
derivation of a Noetherian $\Q$-algebra $R$ and $I$ is an ideal of
$R$ with $\delta(I)\subseteq I$ then $\delta(P)\subseteq P$ for each
prime ideal $P$ of $R$ that is minimal over $I$ and the fact that
$I$ must contain a finite product of such ideals $P$.
\end{rmks}
\begin{defn} Let $A=\C[x,y,z]$ and let $f\in A$. There is a Poisson bracket
$\{-,-\}_f$ on $A$ given by $\{x,y\}_f=\partial f/\partial z$,
$\{y,z\}_f=\partial f/\partial x$ and $\{z,x\}_f=\partial f/\partial
y$ and such that $f\in \PZ(A)$. For more detail, see \cite{DAJNS}.
We call a Poisson bracket on $A$ {\it exact} (determined by $f$) if
it has the form $\{-,-\}_f$ for some $f\in A$.
\end{defn}
\begin{defns}
\label{quantdeform}
If $T$ is a $\C$-algebra with a central non-unit
non-zero-divisor $t$ such that $B:=T/tT$ is commutative then there
is a Poisson bracket $\{-,-\}$ on $B$ such that $\{\ov u,\ov
v\}=\ov{t^{-1}[u,v]}$ for all $\ov u=u+tT$ and $\ov v=v+tT\in B$. In
this situation, we follow \cite[Chapter III.5]{BGl} in referring to
$T$ as a quantization of the Poisson algebra $B$ and we refer to any
$\C$-algebra of the form $T/(t-\lambda)T$, where $\lambda\in \C$ is
such that the central element $t-\lambda$ is a non-unit in $T$, as a
{\it deformation} of $A$.
\end{defns}
\begin{defns}
We shall use the language of modules over Lie algebras rather than
the equivalent language of representations. A module $M$ over a Lie
algebra $\g$ is a $\C$-vector space endowed with a bilinear product
$[-,-]_M:\g\times M\rightarrow M$ such that
\[
[[a,b],m]_M=[a,[b,m]_M]_M-[b,[a,m]_M]_M
\] for
all  $a,b\in \g$ and all $m\in M$. A subspace $V$ of $M$ is a {\it
submodule} of $M$ if $[a,v]_M\in V$ for all $a\in \g$ and all $v\in
V$ and $M$ is {\it simple} if $M\neq 0$ and its only submodules are
$0$ and $M$.
\end{defns}

\begin{defn}
\label{poissonmodule} Let $A$ be a commutative Poisson algebra with
Poisson bracket $\{-,-\}$. An $A$-module $M$ is a {\it Poisson
$A$-module} if there is a bilinear product $\{-,-\}_M:A\times
M\rightarrow M$ such that the following axioms hold for all $a,b\in
A$ and all $m\in M$:
\begin{enumerate}
\item
$\{\{a,b\},m\}_M=\{a,\{b,m\}_M\}_M-\{b,\{a,m\}_M\}_M$;
\item
$\{a,bm\}_M=\{a,b\}m+b\{a,m\}_M$;
\item
$\{ab, m\}_M=a\{b,m\}_M+b\{a,m\}_M$.
\end{enumerate}
\end{defn}

\begin{rmk}
\label{PMremark} There are two definitions of Poisson module in the
literature. The one given above is as in \cite{farkas} and
\cite{oh}. Other authors, for example \cite{BG}, only require a
bilinear product $\{-,-\}_M$ satisfying
Axiom~\ref{poissonmodule}(ii). Such a product is called a weak
Poisson action in \cite{farkas}.

Axiom~\ref{poissonmodule}(i) says that $M$ is a module over $A$ when
$A$ is viewed as a Lie algebra. Axiom~\ref{poissonmodule}(iii)
requires that $\{-,m\}$ is a derivation from $A$ to the $A$-module
$M$ whereas Axiom~\ref{poissonmodule}(ii) is a derivation-like
property for the map $\{a,-\}_M:M\rightarrow M$. To be more precise,
given a derivation $\delta$ of a commutative ring $R$ and an
$R$-module $M$, we can define a $\delta$-{\it derivation} of $M$ to
be a linear map $\theta:M\rightarrow M$ such that
$\theta(am)=\delta(a)m+a\theta(m)$ for all $a\in A$ and $m\in M$.
Axiom~\ref{poissonmodule}(ii) then requires that $\{a,-\}_M$ is a
$\ham (a)$-derivation of $M$.
\end{rmk}

\begin{defns}
Let $A$ be a Poisson algebra. A submodule $N$ of a Poisson
$A$-module is a {\it Poisson submodule} of $M$ if $\{a,n\}_M\in N$
for all $a\in A$ and all $n\in N$. A non-zero Poisson $A$-module $M$
is {\it simple} if its only Poisson submodules are $0$ and $M$ is
{\it Poisson semisimple} if it is a direct sum of simple Poisson
$A$-modules.

Let $M,M^\prime$ be Poisson $A$-modules. A {\it Poisson module
homomorphism} $f:M\rightarrow M^\prime$ is an $A$-module
homomorphism such that $f(\{a,m\}_M)=\{a,f(m)\}_{M^\prime}$  for all
$a\in A$ and $m\in M$. If such a map $f$ is bijective then
$f^{-1}:M^\prime\rightarrow M$ is also a Poisson module homomorphism
and $f$ is a {\it Poisson module isomorphism}.

For the purposes of this paper, we say that a Poisson algebra $A$ is
\textit{$t$-homogeneous}, for a non-negative integer $t$, if, for,
each positive integer $d$, there are, up to isomorphism, precisely
$t$ simple Poisson $A$-modules of dimension $d$. As indicated in the
introduction, we use the term \textit{$t$-homogeneous algebra} for
an algebra with the analogous associative property.
\end{defns}

\begin{rmk}
\label{automodule} If $f:\g\rightarrow \h$ is a homomorphism of Lie
algebras and $M$ is an $\h$-module then there is a $\g$-module
$\leftexp{f}{M}$ such that, as a set, $\leftexp{f}{M}=M$ and
$[a,m]_{\leftexp{f}{M}}=[f(a),m]_M$ for all $a$ and all $m\in M$. We
shall use the same notation for the corresponding construction for
modules over an associative algebra $T$. In particular, if $M$ is a
simple $T$-module and $\theta$ is an automorphism of $T$, then
$\leftexp{\theta}{M}$ is also a simple $T$-module.

If $A$ is a Poisson algebra, $M$ is a Poisson $A$-module and $\pi$
is a Poisson automorphism of $A$ then, combining the above
constructions, we obtain a Poisson $A$-module $\leftexp{\pi}{M}$,
equal to $M$ as a $\C$-vector space, with $a\cdot m=\pi(a)m$ and
$\{a,m\}_{\leftexp{\pi}{M}}=\{\pi(a),m\}_M$ for all $a\in A$ and all
$m\in M$. If $J=\ann_A(M)$ then
$\pi^{-1}(J)=\ann_A(\leftexp{\pi}{M})$ and $\leftexp{\pi}{M}$ is
simple if and only if $M$ is simple.
\end{rmk}

\begin{rmk}\label{IIJmodule} One natural way in which Poisson modules arise is
whenever $I$ and $J$ are Poisson ideals of $A$ with $I\subseteq J$.
The factor $J/I$ is then a Poisson module with
$\{a,j+I\}_{J/I}=\{a,j\}+I$. It is a routine matter to check that
$\{-,-\}_{J/I}$ is well-defined and that the axioms for a Poisson
module are satisfied. By Remark~\ref{Pideals}(i), $J/I$ is also a
Lie algebra. Every Poisson submodule of $J/I$ is a Lie ideal so if
$J/I$ is simple as a Lie algebra then it is simple as a Poisson
module. If $A$ is affine and $J$ is a Poisson maximal ideal, so that
$A=J+\C$, then the converse is also true because every Lie ideal of
$J/I$ is then a Poisson $A$-submodule.

If $I$ and $J$ are any two Poisson ideals of $A$ then $I/IJ$ and
$J/IJ$ are Poisson modules.
\end{rmk}

\begin{notn}
Let $M$ be a Poisson module over a Poisson algebra $A$ and let
$S\subseteq M$. We shall denote the annihilator of $S$ in $A$, in
the module sense, by $\ann_A(S)$ and we shall denote by $\Pann_A(S)$
the set $\{a\in A:\{a,m\}_M=0\mbox{ for all }m\in S\}$.
\end{notn}

\section{Finite-dimensional simple Poisson modules}

\begin{lemma}\label{basicPmod}
Let $A$ be an affine Poisson algebra and let $M$ be a Poisson
$A$-module. Let $J=\ann_A(M)$.

(i) $J$ is a Poisson ideal of $A$.

(ii) If $M$ is simple then $J$ is a prime ideal of $A$.

(iii) If $M$ is finite-dimensional and simple then $J$ is a maximal
ideal of $A$.

(iv) $\C+J^2\subseteq \Pann_A(M)$.
\end{lemma}
\begin{proof}
(i) Let $a\in A$, $j\in J$, and $m\in M$. Then $jm=j\{a,m\}_M=0$ so,
by Axiom~\ref{poissonmodule}(ii),
$
\{a,j\}m=\{a,jm\}-j\{a,m\}_M=0,
$
whence $\{a,j\}\in J$ and $J$ is a Poisson ideal.

(ii) Suppose that $M$ is simple and let $I$ and $K$ be Poisson
ideals of $A$ such that $IK \subseteq J$, that is, $IKM = 0$. If $k
\in K $ and $m \in M $ then, for $a \in A$,
$
\{a,km\}_M = k\{a,m\}_M + \{a,k\}m  \in KM.
$
Hence $KM$ is a Poisson submodule of $M$ so either $KM=0$, in which
case $K\subseteq J$, or $KM=M$, in which case $IM=IKM=0$ and
$I\subseteq J$. By Remark~\ref{Pideals}(iii), $J$ is a prime ideal
of $A$.

(iii) By (ii), $A/J$ is a domain. As it embeds in the
finite-dimensional $\C$-algebra $\Endo_\C(M)$ it must be isomorphic
to $\C$. Hence $J$ is maximal.

(iv) It is clear from  Axiom~\ref{poissonmodule}(iii), with $a=b=1$,
that $\{1,m\}_M=0$ for all $m\in M$ so $\C\subseteq \Pann_A(M)$. Let
$m\in M$ and $j,k\in J$. By Axiom~\ref{poissonmodule}(iii),
$
\{jk,m\}_M=j\{k,m\}_M+k\{j,m\}_M=0,
$
whence $J^2\subseteq \Pann_A(M)$.
\end{proof}

\begin{rmk}
If $A, M$ and $J$ are as in \ref{basicPmod} and $I$ is a Poisson
ideal such that $I\subseteq J^2$ then, by \ref{basicPmod}(iv), $M$
becomes a Poisson $A/I$-module in the obvious way: $(a+I)m=am$ and
$\{a+I,m\}_M=\{a,m\}_M$.
\end{rmk}

\begin{notn}
Let $A$ be an affine Poisson algebra. For a finite-dimensional
simple Poisson $A$-module $M$, let $\widehat{M}$ denote its
isomorphism class. Similarly, for a Lie algebra $\g$ and a
finite-dimensional simple $\g$-module $N$, let $\widehat{N}$ denote
its isomorphism class.
\end{notn}

\begin{theorem}
\label{fdsPmod} Let $A$ be an affine generated Poisson algebra.

(i) Let $M$ be a finite-dimensional simple Poisson $A$-module and
let $J=\ann_A(M)$. There is a simple module $M^*$ for the Lie
algebra $\g(J)$ such that $M^*=M$, as $\C$-vector spaces, and
$[j+J^2,m]_{M^*}=\{j,m\}_M$ for all $j\in J$ and $m\in M$.

(ii)
 Let $J$ be a Poisson maximal ideal of $A$ and let $N$ be a finite-dimensional simple
$\g(J)$-module. There exist a simple Poisson $A$-module $N^\dagger$
and a Lie homomorphism $f:A\rightarrow \g(J)$ such that
$N^\dagger=\leftexp{f}{N}$ as a Lie module over $A$ and
$\ann_A(N^\dagger)=J$.

(iii) For all finite-dimensional simple Poisson modules $M$,
$M^{*\dagger}=M$. For all Poisson maximal ideals $J$ of $A$ and all
finite-dimensional simple $\g(J)$-modules $N$, $N^{\dagger*}=N$.

(iv) The procedures in (i) and (ii) establish a bijection $\Gamma$
from the set of isomorphism classes of finite-dimensional simple
Poisson modules over $A$ to the set of pairs $(J,\widehat{N})$,
where $J$ is a Poisson maximal ideal of $A$ and $N$ is a
finite-dimensional simple $\g(J)$-module, given by
$\Gamma(\widehat{M})=(\ann_A(M),\widehat{M^*})$.
\end{theorem}

\begin{proof}
(i) By Lemma~\ref{basicPmod}(iv), $J^2\subseteq \Pann_A(M)$ and, by
Lemma~\ref{basicPmod}(iii), $J$ is a Poisson maximal ideal of $A$,
whence, $A$ being affine, $A=J\oplus\C$, as vector spaces. Being a
module over the Lie algebra $A$, $M$ is also a module over the Lie
algebra $J$ and, as $J^2\subseteq \Pann_A(M)$, over $\g(J)$. Let
$M^*$ denote this $\g(J)$-module. Thus $[j+J^2,m]_{M^*}=\{j,m\}_M$
for all $j\in J$ and $m\in M$.

Let $V$ be a non-zero submodule of $M^*$. Then $V$ is a Lie
$J$-submodule of $M$. By Lemma~\ref{basicPmod}(iv), $\C\subseteq
\Pann_A(M)$, so, as $A=J\oplus\C$, $V$ is a Lie $A$-submodule of
$M$. As $A=J\oplus\C$ and $JM=0$, $V$ is also an $A$-submodule of
$M$ in the associative sense so it is a non-zero Poisson
$A$-submodule of $M$. Hence $V=M$, as Poisson $A$-modules, and
$V=M^*$, as $\g(J)$-modules. Therefore $M^*$ is a simple
$\g(J)$-module.

 (ii) As $J$ is maximal and $A$ is affine, $A/J\simeq \C$ and
$A=\C\oplus J$ as vector spaces. Let $\pi_J:A\rightarrow J$ and
$\pi_\C:A\rightarrow \C$ be the projections and let
$\phi:J\rightarrow \g(J)$ be the natural Lie homomorphism. Let
$a,b\in A$. Then
\[
\{a,b\}=\{\pi_J(a)+\pi_\C(a),\pi_J(b)+\pi_\C(b)\}=\{\pi_J(a),\pi_J(b)\}.
\]
It follows, since $J$ is a Poisson ideal, that $\{A,A\}\subseteq J$
and that $\pi_J$ is a Lie homomorphism.  Being a $\C$-module, $N$
has the structure of an $A$-module, with $a\cdot n=\pi_\C(a)n$ for
all $a\in A$ and all $n\in N$, and hence with $\ann_A(N)=J$. Let
$N^\dagger$ denote $N$ equipped with this $A$-module structure and
the Lie $A$-module structure of $\leftexp{f}{N}$, where $f$ is the
Lie homomorphism $\phi\pi_J:A\rightarrow \g(J)$. For $a\in A$ and
$n\in N^\dagger$, write
$\{a,n\}_{N^\dagger}=[a,n]_{\leftexp{f}{N}}=[f(a),n]_N$. To check
that $N^\dagger$ is a Poisson $A$-module, it remains to check
Axioms~\ref{poissonmodule}(ii) and (iii).

Let $a,b\in A$ and $n\in N$ and let $\alpha=\pi_\C(a),
\beta=\pi_\C(b), j=\pi_J(a)$ and $k=\pi_J(b)$. For
Axiom~\ref{poissonmodule}(ii),
\[
\{a,b\cdot n\}_{N^\dagger}=[f(a),\beta n]_{N}= \beta[f(a),n]_N=
b\cdot\{a,n\}_{N^\dagger}
\]
and, as $\{a,b\}\in J$, $\{a,b\}\cdot n=0$. Therefore
Axiom~\ref{poissonmodule}(ii) holds.

For Axiom~\ref{poissonmodule}(iii),  note that $ab=jk+\alpha k+\beta
j+\alpha\beta$, $f(a)=f(j)$, $f(b)=f(k)$ and
$f(jk)=0=f(\alpha\beta)$. Hence
\begin{align*}
\{ab,n\}_{N^\dagger}=&[f(ab),n]_N\\
=&[f(\alpha k+\beta j),n]_N\\
=&\alpha [f(k),n]_N+\beta [f(j),n]_N\\
=&\alpha[f(b),n]_N+\beta [f(a),n]_N\\
=&a\cdot\{b,n\}_{N^\dagger}+b\cdot\{a,n\}_{N^\dagger}.
\end{align*}
Thus $N^\dagger$ is a Poisson $A$-module with $\ann_A(N^\dagger)=J$.

To see that $N^\dagger$ is simple as a Poisson $A$-module, let $W$
be a non-zero Poisson submodule of $N^\dagger$ and let $w\in W$,
$\phi(j)\in \g(J)$. Then
$[\phi(j),w]_N=[f(j),w]_N=\{j,w\}_{N^\dagger}\in W$. Thus $W$ is a
$\g(J)$-submodule of $N$, and, by the simplicity of $N$,
$W=N^\dagger$.

(iii) For $m\in M$ and $a\in A$,
\[
\{a,m\}_{M^{*\dagger}}=[\phi\pi_J(a),m]_{M^*}=\{\pi_J(a),m\}_M=\{a,m\}_M,
\]
where the last step uses the fact that $\C\subseteq \Pann_A(M)$.
Also
 $a\cdot m=\pi_\C(a)m=am$. Therefore $M^{*\dagger}=M$.

If $j\in J$ and $n\in N$,
\[
[\phi(j),n]_{N^{\dagger*}}=\{j,n\}_{N^\dagger}=[f(j),n]_N=[\phi(j),n]_N.
\] Thus $N^{\dagger*}=N$.

(iv) If $g:M\rightarrow M_1$ is an isomorphism of finite-dimensional
simple Poisson modules then it is clear that $J=\ann_A(M_1)$ and
that $g:M^*\rightarrow M_1^*$ is an isomorphism of Lie
$\g(J)$-modules. Also, if $h:N\rightarrow N_1$ is an isomorphism of
finite-dimensional simple Lie $\g(J)$-modules, for some Poisson
maximal ideal $J$, then $h:N^\dagger\rightarrow N_1^\dagger$ is an
isomorphism of Poisson $A$-modules. The result now follows from
(iii).
\end{proof}

\begin{example}
Let $\g$ be a finite-dimensional complex Lie algebra with basis
$x_i$, $1\leq i\leq n$, and let $A:=\C[x_1,x_2,\ldots,x_n]$. The
Kirillov-Kostant Poisson bracket \cite[III.5.5]{BGl} on $A$ is such
that $\{x_i,x_j\}=[x_i,x_j]$. If $[\g,\g]=\g$, in particular if $\g$
is semisimple, then the only Poisson maximal ideal is
$J=x_1A+x_2A+\ldots+x_nA$ and $\g(J)\simeq \g$. Thus the
finite-dimensional simple Poisson $A$-modules correspond exactly to
the finite-dimensional simple $\g$-modules.

At the other extreme, if $\g$ is abelian then the Poisson bracket on
$A$ is zero, every maximal ideal of $A$ is Poisson and each
finite-dimensional simple Poisson $A$-module $M$ is one-dimensional,
with basis $\{m\}$ and a pair of $n$-tuples $\alpha_i$ and $\beta_i$
such that, for $1\leq i\leq n$, $x_im=\alpha_i m$ and
$\{x_i,m\}_M=\beta_i m$. For $f\in A$,
\[\{f,m\}_M=\sum_{i=1}^n \beta_i \frac{\partial f}{\partial
x_i}(\alpha_1,\alpha_2,\ldots,\alpha_n)m.\]
\end{example}

\section {\bf Exact examples.}
\subsection{Some generalities} In this section we consider Poisson algebras of the form
$A=\C[x,y,z]$, with an exact bracket $\{-,-\}_f$, $f\in A$, and
their factors. For $\lambda\in \C$,  $A/(f-\lambda)A$ will be
written $A_\lambda$ and we abuse notation by writing the images of
$x, y$ and $z$ in $A_\lambda$ as $x,y,z$ respectively, rather than
$\ov x,\ov y,\ov z$. The surface $\{(a,b,c)\in
\C^3:f(a,b,c)=\lambda\}$, of which $A_\lambda$ is the coordinate
ring, will be written $S_\lambda$. If $J$ is a Poisson maximal ideal
of $A$ containing $f-\lambda$ then the Poisson maximal ideal
$J/(f-\lambda)A$ of $A_\lambda$ will be written $\J$. In several
cases, $A_0$ is the Poisson algebra of invariants for the action of
finite group of Poisson automorphisms of a simple Poisson algebra
$B$ of Krull dimension $2$ and has a deformation that is the algebra
of invariants for the action of a finite group of automorphisms of a
simple noncommutative deformation $R$ of $B$.

Let $J=(x-\alpha,y-\beta,z-\gamma)$ be a maximal ideal of $A$ and
let $\lambda=f(\alpha,\beta,\gamma)$. Thus $f-\lambda\in J$ and
determines the same exact bracket on $A$ as $f$. The maximal ideal
$J$ is a Poisson ideal, for the exact bracket $\{-,-\}_f$, if and
only if $\partial f/\partial x\in J$, $\partial f/\partial y\in J$
and $\partial f/\partial z\in J$ if and only if $f-\lambda\in J^2$
if and only if $(\alpha,\beta,\gamma)$ is a singularity of the
surface $S_\lambda$.

Now suppose that $J$ is Poisson and let $\lambda$ be as before. The
maximal Poisson ideals $J$, of $A$, and $\J$ of $A_\lambda$,
determine the Lie algebras $\g(J)$ and  $\g(\ov J)$.  As
$f-\lambda\in J^2$, these Lie algebras are isomorphic so we shall
use the single notation $\g(J)$ to denote them. When
Theorem~\ref{fdsPmod} is applied to determine the finite-dimensional
simple Poisson $A$-modules, the finite-dimensional simple Poisson
$A_\lambda$-modules will automatically be determined, being the
finite-dimensional simple Poisson $A$-modules annihilated by those
Poisson maximal ideals that contain $f-\lambda$.

If $S_\lambda$ is not smooth let $\sing(S_\lambda)$ denote its
singular locus. The partition of $\C^3$ into the surfaces
$S_\lambda$ can be refined to partition $\C^3$ into (i) the smooth
surfaces $S_\lambda$, (ii) the punctured surfaces
$S_\lambda\backslash \sing(S_\lambda)$ and (iii) the points that are
singularities for some $S_\lambda$. These are the so-called
symplectic leaves of $A$, see \cite[p. 274]{BGl} or \cite{BG}.

The Lie algebras $\g(J)$ that arise in many of the examples that
follow are $3$-dimensional and several of them are isomorphic to
$sl_2$. We shall not give explicit details of such isomorphisms but
will appeal, implicitly, to the fact \cite[3.2.4]{erdmann} that, up
to isomorphism, $sl_2$ is the unique $3$-dimensional Lie algebra
$\s$ over $\C$ with $[\s,\s]=\s$.

In the next three subsections, we consider exact Poisson brackets on
$A$ associated to Poisson subalgebras of invariants of the simple
Poisson algebras $\C[x_1,x_2]$, with $\{x_1,x_2\}=1$, $\C[x_1^{\pm
1},x_2^{\pm 1}]$, with $\{x_1,x_2\}=x_1x_2$, and $\C[x_1^{\pm
1},x_2]$, with $\{x_1,x_2\}=x_1$.

\subsection{The Kleinian singularity of type $A_1$ }
\label{KleinianA1} Let $f=z^2-xy$, so that  $\{x,y\}=2z$, $\{
y,z\}=-y$ and $\{z,x\}=-x$. The Poisson bracket $\{-,-\}_f$ is
equivalent to the Kirillov-Kostant Poisson bracket on $A$ for the
Lie algebra $sl_2$. The coordinate ring $A_0=A/fA$  has a unique
singularity at $0$, the Kleinian singularity of type $A_1$, and has
a well-known interpretation as a subalgebra of invariants. Let $B$
be the Poisson algebra $B=\C[x_1,x_2]$, with $\{x_1,x_2\}=1$, and
let $\pi$ be the Poisson automorphism of $B$ such that
$\pi(x_1)=-x_1$, $\pi(x_2)=-x_2$. The Poisson algebra of invariants
$B^\pi$ may be identified with $A_0$, with $x=x_1^2/2$, $y=x_2^2/2$
and $z=x_1x_2/2$, so that $xy=z^2$, $ \{x,y\}=2z$, $\{y,z\}=-y$ and
$\{z,x\}=-x$.

There is a unique Poisson maximal ideal, $J:=(x,y,z)$ and $\g(J)$ is
isomorphic to $sl_2$. Hence, by Theorem~\ref{fdsPmod}, the Poisson
algebras $A$ and $A_0$ are both $1$-homogeneous.

\begin{prop}\label{Aquant}
Let $U$ denote the $\C$-algebra generated by $x, y, z$ and $t$
subject to $t$ being central and the further relations
\begin{equation}
xy=yx+2tz+\tfrac{1}{2}t^2,\quad yz=zy-ty,\quad zx=xz-tx.\label{xyyx} \\
\end{equation}

(i) The algebra $U$ is a quantization of $A$, with the Poisson
bracket $\{-,-\}_f$, and $U_1:=U/(t-1)U$ is a deformation of $A$.

(ii) The elements
$g:=z^2-\tfrac{1}{2}tz-yx=z^2+\frac{3}{2}zt+\frac{1}{2}t^2-xy$ and
$p:=z^2-\tfrac{1}{2}z-yx=z^2+\frac{3}{2}z+\frac{1}{2}-xy$ are
central in $U$ and $U_1$ respectively.

(iii) $U_1$ is $1$-homogeneous and, for each $d\geq 1$, the unique
simple module of dimension $d$ is annihilated by $p-\sigma_d$, where
$\sigma_d:=\frac{4d-1}{16}$.
\end{prop}
\begin{proof}
(i) is a routine calculation. Both $U$ and $U_1$ are examples of the
iterated skew polynomial rings introduced in \cite{itskew} and
subsequently generalized and named ambiskew polynomial rings
\cite{jorbav,ambi}. In the notation of \cite{itskew}, take
$A=\C[z,t]$, $\alpha(t)=t$, $\alpha(t)=t-1$ and
$u=z^2+\frac{3}{2}zt+\frac{1}{2}t^2$ to obtain $U$ and make the
obvious modifications, setting $t=1$, to obtain $U_1$. (ii) is then
a consequence of \cite[1.2]{itskew} and (iii) is obtained by
applying \cite[3.11 and 3.6]{itskew} or \cite[2.3 and 2.1]{fds}.
\end{proof}

By Proposition \ref{Aquant}, there is, up to isomorphism, a
 bijection, preserving dimension, between the finite-dimensional simple
modules for $A$ and the finite-dimensional simple modules of the
deformation $U_1$, albeit without explicit machinery for passing
between the two. A similar situation exists for several of the
examples of exact brackets that follow. It would be interesting to
know for which Poisson brackets such a bijection exists. To see that
there is no general result of this form, for Poisson algebras and
their deformations, we consider the Poisson algebra $A_0=B^\pi$.

The algebra $W=\C[x_1,x_2,t:x_1x_2-x_2x_1=t]$, with $t$ central, is
a quantization of the Poisson algebra $B$ and the first Weyl algebra
$W_1:=W/(t-1)W$ is an associated deformation. The $\C$-automorphisms
$\theta$ of $W$ and $\theta_1$ of $W_1$ such that $x\mapsto-x$,
$y\mapsto -y$ and, for $\theta$, $t\mapsto t$ are a quantization and
deformation of $\pi$ in the sense of \cite{DAJNS}. The algebras
$W^\theta$ and $W_1^{\theta_1}$ belong to families of quantizations
and deformations of $B^\pi=A_0$.

\begin{prop}\label{T}
Let $\sigma\in \C$. With $U,U_1,g,p$ and $\sigma_d$ as in
Proposition \ref{Aquant} and $W, W_1, \theta$ and $\theta_1$ as
above, let $T_\sigma=U/(g-\sigma t)U$ and
$T_{\sigma,1}=U_1/(p-\sigma)U_1$. Thus $T_\sigma$ is generated by
$x, y, z$ and a central generator $t$ subject to \eqref{xyyx} and
either of the relations
\[
yx=z^2-\tfrac{1}{2}tz-\sigma t\label{yx},\quad
xy=z^2+\tfrac{3}{2}zt+\tfrac{1}{2}t^2-\sigma t
\] and the relations for $T_{\sigma,1}$ are obtained from
those for $T_{\sigma}$ by setting $t=1$.

(i) $T_{\sigma}$ and $T_{\sigma,1}$ are, respectively, a
quantization and deformation of $A_0$.

(ii) $T_{0}=W^\theta$ and $T_{0,1}=W_1^{\theta_1}$.

(iii) $T_{\sigma,1}$ is simple if and only if $\sigma\neq \sigma_d$
for all $d\geq 1$.

(iv) If $\sigma=\sigma_d$ for some, necessarily unique, $d$ then
$T_{\sigma,1}$ has a unique finite-dimensional simple module $M$ and
$\dim M=d$.
\end{prop}
\begin{proof}
(i) is a straightforward calculation. (ii) is well-known for $W_1$,
for example see \cite{hodges}, and the result can be lifted to $W$.
Alternatively it is easy to check that, in $W$, $x:=x_1^2/2$,
$y:=x_2^2/2$, $z:=x_1x_2/2$ and $t$ are invariants for $\theta$,
that they generate $W^\theta$, that they satisfy the defining
relations for $T_0$ and, using standard methods of Gelfand-Kirillov
dimension, that if $\chi:T_0\rightarrow W^\theta$ is the resulting
surjective algebra homomorphism then $\ker \chi=0$. The same
approach can be used for $W_1$. (iii) is a routine application of
\cite[Theorem 6.1]{prim}. In the notation of \cite[5.1]{prim},
$A=\C[z]$, $\alpha(z)=z-1$ and
$T_{\sigma,1}=T(z^2+\frac{3}{2}z+\frac{1}{2}-\sigma)$. Finally (iv)
is a re-statement of Proposition~\ref{Aquant}(iii).
\end{proof}

As has been observed by Alev and Farkas in \cite{alevfarkas}, the
existence of the Poisson ideal $\J$ of $A_0$ shows that, in contrast
to the situation for the deformation $W_1^{\theta_1}$, which is
simple as a consequence of \cite[Theorem 28.3(ii)]{passman}, Poisson
simplicity is not preserved by the passage from $B$ to $B^\pi\simeq
A_0$. Although $B^\pi\simeq A_0$ is $1$-homogeneous, its deformation
$W^\theta\simeq T_{0,1}$, has no finite-dimensional simple modules.
However, we see from Proposition~\ref{T}, that, for each $d\geq 1$,
there exists a unique $\sigma=\sigma_d$ such that the deformation
$T_{\sigma,1}$ of $A_0$ has a $d$-dimensional simple module and that
there is a  bijection, preserving dimension, between the isomorphism
classes of finite-dimensional simple Poisson modules for $A_0$ and
the union of the isomorphism classes of finite-dimensional simple
modules over the deformations $T_{\sigma,1}$, taken over all
$\sigma\in \C$.

\subsection{Invariants of the $2$-torus}
\label{torus} Let $f=xyz-x^2-y^2-z^2+4$ and let $B$ be the
coordinate ring $\C[x_1^{\pm 1},x_2^{\pm 1}]$ of the torus, endowed
with the Poisson bracket such that $\{x_1,x_2\}=x_1x_2$. If $\pi$ is
the $\C$-automorphism of $B$ such that $\pi(x_1)=x_1^{-1}$ and
$\pi(x_2)=x_2^{-1}$ then $\pi$ is Poisson and $B^\pi$ is generated
by $x:=x_1+x_1^{-1}, y:=x_2+x_2^{-1}$ and
$z:=x_1x_2+x_1^{-1}x_2^{-1}$ subject to the single relation
$xyz-x^2-y^2-z^2+4=0$. It is known that $B^\pi\simeq A_0$. For
example, see \cite[Example 3.5]{lor}, although there the base ring
is $\Z$ rather than $\C$. In the Poisson algebra $B^\pi$,
\[\{x,y\}=xy-2z,\; \{y,z\}=yz-2x\mbox{ and }\{z,x\}=zx-2y.\] Thus
there is a Poisson isomorphism $B^\pi\simeq A_0$, where the Poisson
bracket on $A_0$ is induced by the exact bracket $\{-,-\}_f$ on $A$.

There are five Poisson maximal ideals in $A$:
\begin{gather*}
J_1:=xA+yA+zA,\\
J_2:=(x-2)A+(y-2)A+(z-2)A,\\
J_3:=(x-2)A+(y+2)A+(z+2)A,\\
J_4:=(x+2)A+(y-2)A+(z+2)A,\\ J_5:=(x+2)A+(y+2)A+(z-2)A.
\end{gather*}
All but two of the surfaces $S_\lambda$ are smooth. The
singularities corresponding to $J_2, J_3, J_4, J_5$ lie on $S_0$ and
$S_4$ has a singularity at $(0,0,0)$, corresponding to $J_1$.

The Lie algebra $\g({J_1})$ has basis (the images of) $x, y, z$ and
$[y,x]=2z$, $[z,y]=2x$, $[x,z]=2y$. Consequently $\g(J_1)\simeq
sl_2$. It follows that, for each $d\geq 1$, $A$ has a unique
$d$-dimensional simple Poisson module annihilated by $J_1$ and that
$A_4$ is $1$-homogeneous.

For each of the generators $g=x, y$ or $z$ of $A$, there is Poisson
automorphism $\theta_g$ of $A$ fixing $g$ and with the other two
generators eigenvectors for $-1$. As $J_3=\theta_x(J_2)$,
$J_4=\theta_y(J_2)$ and $J_5=\theta_z(J_2)$, the four Lie algebras
$\g({J_i})$, $2\leq i\leq 5$, are isomorphic.

Let $u=x-2$, $v=y-2$, $w=z-2$ and $J=J_2=u A +v A+w A$. Then, in
$J$, $ \{u,v\}=u v+2u+2v- 2w$, $\{v,w\}=v w-2u+2v+ 2w$ and
$\{w,u\}=u w+2u-2v+2w$, so, in $\g({J})$,
\[
[u,v]=2u+2v-2w,\;[v,w]=-2u+2v+ 2w,\;[w,u]=2u-2v+2w,
\]
from which it follows that $\g({J})\simeq sl_2$. Therefore, for each
$d\geq 1$, $A$ has a unique $d$-dimensional simple Poisson module
annihilated by $J_2$. The same conclusion holds for $J_3, J_4$ and
$J_5$, for each of which the $d$-dimensional simple Poisson module
annihilated by $J_i$ is obtained from that annihilated by $J_2$ by
the construction described in \ref{automodule}. Thus $A$ is
$5$-homogeneous and $A_0$ is  $4$-homogeneous.

For $q\in \C$, the $\F$-algebra $T_q$ generated by $x, y$ and $z$
subject to the relations
\begin{align}
xy-qyx&=(1-q^2)z,\label{qx12} \\
yz-qzy&=(q^{-1}-q)x,\label{qx23} \\
zx-qxz&=(q^{-1}-q)y.\label{qx31}  \end{align} is a deformation of
$A$, under $\{-,-\}_f$ and there is a central element
 $p=zyx-qz^2-q^{-2}y^2-x^2+2(1+q^{-2})$ of $T_q$ such that
 $T_q/pT_q$ is a deformation of $A_0$. See \cite{DAJNS} for details.
 When $q\neq 1$, the algebra $T_q$ is isomorphic to the algebra
 denoted
 $U_q^\prime(so_3)$ in \cite{mp2}. It is shown in \cite{mp2}, and,
 by an independent method in \cite{nongthesis},
  that, provided $q$ is not a root of unity,
$U_q^\prime(so_3)$ is $5$-homogeneous, reflecting the situation for
the  Poisson algebra $A$.

Using the details of the modules in \cite{mp2} or \cite{nongthesis},
one can check that, for each $d\geq 1$, there exist $\lambda_d,
\mu_d\in \C$, all distinct, such that $p-\lambda_d$ annihilates four
of the $d$-dimensional simple $T_q$-modules and $p-\mu_d$
annihilates the fifth. If $M$ and $N$ are any two of the four that
are annihilated by $p-\lambda_d$ then there exists a
$\C$-automorphism $\theta$ of $T_q$ such that $\theta(p)=p$ and
$N=\leftexp{\theta}{M}$. Each of the algebras $T_q/(p-\lambda)T_q$,
$\lambda\in \C$, is a deformation of $A_\alpha$ for some $\alpha\in
\C$. If $\lambda=\lambda_d$ for any $d$ then $\alpha=0$ and if
$\lambda=\mu_d$ for any $d$ then $\alpha=4$. Thus, as in
Subsection~\ref{KleinianA1}, for each of $A_0$ and $A_4$, there is a
 bijection, preserving dimension, between the isomorphism classes of
finite-dimensional simple Poisson modules and the union, over a
family of deformations, of the isomorphism classes of
finite-dimensional simple modules.

With $A_0$ as above, let $R$ be the group algebra $\C V$ of the
Klein four-group $V=\{1,a,b,c\}$, considered as a Poisson algebra
with the zero bracket. There is a surjective Poisson homomorphism
$\psi:A_0\rightarrow R$, given by $x\mapsto 2a, y\mapsto 2b$ and
$z\mapsto 2c$. Every $R$-module $M$ is a Poisson module, with
$\{-,-\}_M=0$, and the finite-dimensional simple (Poisson)
$R$-modules are the four one-dimensional modules $R/\psi(\J_i)$,
$2\leq i\leq 5$. For each $i$, $\psi(\J_i)$ is idempotent so
$\g({\psi(\J_i)})=0$. This conforms to Theorem~\ref{fdsPmod}, when a
simple module $M$ over the zero Lie algebra is interpreted as a
one-dimensional vector space $\C v$ with $[0,v]_M=0$.

\subsection{Invariants of $\C[x_1^{\pm 1},x_2]$} \label{Plocenv} Let
$f=x(4-z^2)+y^2$ and let $D$ be the Poisson algebra $\C[x_1,x_2]$
with the Poisson bracket such that $\{x_1,x_2\}=x_1$. This is the
Kirillov-Kostant bracket for the $2$-dimensional non-abelian
solvable Lie algebra over $\C$ and it is a routine matter to check
that every non-zero Poisson ideal contains a power of $x_1$. The
Poisson bracket extends to the localization $B=\C[x_1^{\pm 1},x_2]$
which is a simple Poisson algebra and has a Poisson automorphism
$\pi$ such that $\pi(x_1)=x_1^{-1}$ and $\pi(x_2)=-x_2$. The Poisson
algebra of invariants \cite[Example 3.2 and 3.6]{DAJNS} is generated
by $x:=x_2^2$, $y:=x_2(x_1-x_1^{-1})$ and $z:=x_1+x_1^{-1}$, subject
to the relation $x(4-z^2)+y^2=0$, and the Poisson bracket is given
by
\[\{x,y\}=-2xz,\, \{y,z\}=4-z^2\mbox{ and
}\{z,x\}=2y.
\] Thus there is an isomorphism of Poisson algebras between $B^\pi$ and $A_0$
where the Poisson bracket on $A_0$ is induced by the exact bracket
$\{-,-\}_f$ on $A$.

There are two Poisson maximal ideals in $A$, $J_1:=xA+yA+(z-2)A$ and
$J_2:=xA+yA+(z+2)A$ and $f\in J_1\cap J_2$. If $\phi$ is the Poisson
automorphism  of $A$ such that $\phi(x)=x$, $\phi(y)=-y$ and
$\phi(z)=-z$ then $\phi(J_1)=J_2$ so $\g({J_1})\simeq \g({J_2})$.
Writing $w$ for $z-2$, the Lie algebra $\g({J_1})$ is
$3$-dimensional, with a basis consisting of (the images of) $x, y$
and $w$ and with
\[[x,y]=-4x,\; [y,w]=-4w\mbox{ and }[w,x]=2y.\] Hence $\g({J_1})\simeq
sl_2$ and, by Theorem~\ref{fdsPmod}, $A$ and $A_0$ are
$2$-homogeneous.

For the Poisson bracket $\{-,-\}_f$, $A$ has a quantization $T$
generated by $x, y, z$ and $t$ subject to $t$ being central and the
additional relations
\begin{gather*}
xy=yx-2tzx+3t^2y+2t^3z,\\
yz=zy-tz^2+4t,\quad xz=zx-t^2z-2ty,
\end{gather*}
and with a central element $g=(4-z^2)x+y^2+3tzy+t^2z^2+4t^2$.
Replacing $t$ by $1$, we obtain a deformation $T_1$ with central
element $p=(4-z^2)x+y^2+3zy+z^2+4.$ The factors $T/gT$ and
$T_1/pT_1$ are, respectively, a quantization and deformation of
$A_0$. Both $T$ and $T_1$ are iterated skew polynomial rings over
$\C$ and $T_1/pT_1$ is simple. For more detail, see \cite{DAJNS}. As
far as we know, the classification of the finite-dimensional simple
$T_1$-modules is not in the literature. The classification of the
finite-dimensional simple Poisson modules for $A$ suggests that it
should be $2$-homogeneous. This is indeed the case. The argument
involves passing to the localization at the powers of $z^2-4$,
replacing the generator $x$ by $(z^2-4)^{-1}p$ and observing that
the localization is a polynomial ring over a simple ring.
Consequently, every finite-dimensional simple module is generated by
an eigenvector for the action of $z$ with eigenvalue $\pm 2$. We
leave the technical details to the interested reader but point out
that, for one-dimensional modules, the calculation is
straightforward, there being two maximal ideals of codimension $1$,
generated by $x+\frac{1}{4}$, $y+1$ and $z-2$ and by
$x+\frac{1}{4}$, $y-1$ and $z+2$.

The correspondence between finite-dimensional simple Poisson
$A_0$-modules and finite-dimensional simple modules over
deformations is similar to that in Subsections~\ref{KleinianA1} and
\ref{torus}. For example, each of the one-dimensional $T_1$-modules
is annihilated by $p-3$ so the two one-dimensional simple Poisson
$A_0$-modules correspond to one-dimensional simple modules over the
deformation $T_1/(p-3)T_1$ of $A_0$,
 rather than the simple deformation $T_1/pT_1$.

\subsection{The quantized enveloping algebra $U_q(sl_2)$}
Here we consider a Poisson algebra associated to the
quantized enveloping algebra $U_q(sl_2)$. This requires two
adjustments to the notion of an exact Poisson bracket on
$A=\C[x,y,z]$. Firstly, we may replace $A$ by any overring
$A^\prime$, such as a localization, to which the three partial
derivatives extend.  Secondly, given a Poisson bracket $\{-,-\}$ on
$A$ (or $A^\prime$), and an element $a\in A$ (or $A^\prime$), it is
a routine exercise on the Jacobi identity to verify that $a\{-,-\}$
is a Poisson bracket $\{-,-\}$ on $A$ (or $A^\prime$). This property
of $A$ does not extend to higher dimensions. Here we take
$A^\prime=\C[x,y,z^{\pm 1}]$, $f=xy+z+z^{-1}$ and $a=2z$ to obtain
Poisson bracket $\{-,-\}=2z\{-,-\}_f$ on $A^\prime$, for which
\begin{equation}
\label{Pbuqsl} \{x,y\}=2(z-z^{-1}),\quad \{y,z\}=2zy\;\mbox{ and
}\{z,x\}=2zx.\end{equation}
 There are two Poisson maximal ideals in $A^\prime$,
$J_1:=xA^\prime+yA^\prime+(z-1)A^\prime$ and
$J_2:=xA^\prime+yA^\prime+(z+1)A^\prime$, and they are transposed by
the Poisson automorphism $\phi$ of $A^\prime$ such that $\phi(x)=x$,
$\phi(y)=-y$ and $\phi(z)=-z$. Hence $\g({J_2})\simeq \g({J_1})$.
Setting $w=z-1$, the Lie algebra $\g({J_1})$ is generated by (the
images of) $x$, $y$ and $w$ and
\[
[x,y]=2w,\; [y,w]=-2y\mbox{ and }[w,x]=-2x,\] whence
$\g({J_1})\simeq sl_2$. By Theorem~\ref{fdsPmod}, $A^\prime$ is
$2$-homogeneous.

For $q\in \C^*$, let $V_q$ be the $\C$-algebra generated by $x, y$
and $z^{\pm 1}$ subject to the relations
\[xy-yx=(q-q^{-1})({z-z^{-1}}),\quad yz=q^2zy;\quad xz=q^{-2}zx.
\]
Then $V_q$ is a deformation of $A^\prime$, the corresponding
quantization being obtained by regarding $q$ as a central invertible
indeterminate. If $q^2\neq 1$, there is an isomorphism $\psi$ from
$V_q$ to $U_q(sl_2)$, as presented in \cite{BGl,Kassel}, with
$\psi(z)=K$, $\psi(x)=(q-q^{-1})E$ and $\psi(y)=(q-q^{-1})F$.  The
$2$-homogeneity of the Poisson algebra $A^\prime$ reflects the
$2$-homogeneity,  when $q$ is not a root of unity, of $U_q(sl_2)$,
details of which may be found in \cite[VI.3]{Kassel} or
\cite[Chapter I.4]{BGl}.

For $\lambda\in \C$, let
$A^\prime_\lambda=A^\prime/(f-\lambda)A^\prime$. Then $f-2\in J_1$
and $f+2\in J_2$, so $A^\prime_2$ and $A^\prime_{-2}$ are both
$1$-homogeneous. As in earlier examples, to reflect this in the
representation theory of the deformation $V_q$, when $q$ is not a
root of unity, one has to take the union of the isomorphism classes
of simple modules over the factors $V_q/(c-\lambda)V_q$, where
$\lambda\in \C$ and $c$ is the image, under $\psi^{-1}$, of the
quantized Casimir element of $U_q(sl_2)$.

Ito, Terwilliger and Weng \cite{teretal} have shown that $U_q(sl_2)$
is isomorphic to the $\C$-algebra $Y_q$ generated by $x, y, z^{\pm 1
}$ subject to the relations
\[
{q^2xy-yx}={q^2-1},\, {q^2yz-zy}={q^2-1}\mbox{ and
}{q^2zx-xz}={q^2-1}.
\]
This is a deformation of $A^\prime$ with the exact Poisson bracket
$\{-,-\}_g$, where $g=2(x+y+z-xyz)\in A$, for which
\[
\{x,y\}=2(1-xy), \; \{y,z\}=2(1-yz)\mbox{ and } \{z,x\}=2(1-xz).
\]
Here $J_1:=(x-1)A^\prime+(y-1)A^\prime+(z-1)A^\prime$ and
$J_2:=(x+1)A^\prime+(y+1)A^\prime+(t+1)A^\prime$ are the only
Poisson maximal ideals and $\g({J_1})\simeq \g({J_2})\simeq sl_2$.
It can be checked that the isomorphism constructed in \cite{teretal}
between $U_q(sl_2)$ and $Y_q$ determines a Poisson isomorphism
$\eta$ from $A^\prime$, under $\{-,-\}_g$, to $A^\prime$, under
$2z\{-,-\}_f$, with $\eta(x)=1-zy$, $\eta(y)=x-z^{-1}$, $\eta(z)=z$,
$\eta^{-1}(x)=y+z^{-1}$ and $\eta^{-1}(y)=z^{-1}(1-x)$.

There is an alternative $4$-homogeneous version of $U_q(sl_2)$,
appearing, for example, in \cite{mp1,mp2,itskew,smith} and obtained
by replacing the relation $xy-yx=(q-q^{-1})({z-z^{-1}})$ by
$xy-yx=(q^2-q^{-2})({z^2-z^{-2}})$. Again this is a deformation of
$A^\prime$ with an exact bracket similar to that in \eqref{Pbuqsl},
but with $\{x,y\}=4(z^2-z^{-2})$. There are four Poisson maximal
ideals $J$, each such that $\g({J_1})\simeq sl_2$, and $A^\prime$ is
$4$-homogeneous.

\subsection{One-dimensional simple Poisson modules}
In this subsection, we consider some Poisson algebras for which all
the finite-dimensional simple Poisson modules are one-dimensional.
\begin{example}
Here we consider an example where the singularities of $f$ are not
isolated. Let $f=xy^2-z^2$, so that $A_0$ is the coordinate ring of
the Whitney umbrella. Then $\{x,y\}_f=-2z$, $\{y,z\}_f=y^2$ and
$\{z,x\}_f=2xy$. The only value of $\lambda$ for which $S_\lambda$
is not smooth is $0$, the singularities of $S_0$ are at the points
$(\alpha,0,0)$, $\alpha\in\C$ and the Poisson maximal ideals have
the form $J_\alpha=uA+yA+zA$, where $u=x-\alpha$. The Lie algebra
$\g({J_\alpha})$ has basis (the images of) $u,y,z$ with $[u,y]=-2z$,
$[y,z]=0$ and $[z,u]=2\alpha y.$ This algebra is solvable and if
$\alpha=0$ it is isomorphic to the $3$-dimensional Heisenberg Lie
algebra. By \cite[Corollary 1.3.13]{dix}, every finite-dimensional
simple $\g({J_\alpha})$-module $N$ is one-dimensional and
annihilated by $[\g({J_\alpha}),\g({J_\alpha})]$. Therefore every
finite-dimensional simple Poisson $A$-module $M$ is one-dimensional,
$M=\C v$, say, with $xv=\alpha v, yv=zv=0, \{x,v\}_M=\mu v,
\{y,v\}_M=\rho v$ and $\{z,v\}_M=0$, for some $\alpha,\mu,\rho\in
\C$ with $\alpha\rho=0$.
\end{example}

\begin{example} \label{KleinianA3up} Let $n\geq 2$ and let $f=f_n=z^n-xy$, so
that $A_0$ is the coordinate ring for the Kleinian singularity of
type $A_{n-1}$. We shall see that the situation for $n>2$ is quite
different to that in Subsection~\ref{KleinianA1} where $n=2$.

Let $B=\C[x_1,x_2]$, with the Poisson bracket $\{x_1,x_2\}=1$, and
let $\pi_n$ be the Poisson automorphism of $B$ such that
$\pi_n(x_1)=\varepsilon_n x_1$ and $\pi_n(x_2)=\varepsilon_n^{-1}
x_2$, where $\varepsilon_n$ is a primitive $n$th root of unity. For
purposes of normalization, let $a\in\C$ be such that $a^2=n^{n}$.
Then $B^{\pi_n}$ is generated by $x:=x_1^n/a$, $y:=x_2^n/a$ and
$z:=x_1x_2/n$ subject to the relation $xy=z^n$ and $B^{\pi_n}\simeq
A_0$ as $\C$-algebras. In the Poisson algebra $B^{\pi_n}$,
\[\{x,y\}=nz^{n-1},\; \{y,z\}=-y \mbox{ and }\{z,x\}=-x\] so
$B^{\pi_n}\simeq A_0$ as Poisson algebras. The maximal ideal
$J:=xA+yA+zA$ is the unique Poisson maximal ideal of $A$ for
$\{-,-\}_f$.

Now suppose that $n>2$. Then, in contrast to type $A_1$, $z^{n-1}\in
J^2$ so, in $\g(J)$, $[x,y]=0$, $[y,z]=-ny$ and $[z,x]=-nx$. This
algebra is solvable so the finite-dimensional simple Poisson
$A$-modules are one-dimensional. They have the form $M=\C v$ with
$xv=yv=zv=0=\{x,v\}_M=\{y,v\}_M$ and $\{z,v\}_M=\tau v$, $\tau\in
\C$. The same is true for $A_0$.

When $n>2$ is even, $B^{\pi_n}=(B^{\pi_2})^\theta$ for a Poisson
automorphism $\theta$ of $B^{\pi_2}$ and it exhibits some surprising
behaviour of simple Poisson modules on the passage from a Poisson
algebra to a Poisson subalgebra of invariants for the action of a
finite group of Poisson automorphisms.

Recall, from Subsection~\ref{KleinianA1}, that $B^{\pi_2}$ is
$1$-homogeneous. Take $n=4$ and rewrite $x,y$ and $z$ as $u,v$ and
$w$, respectively, to distinguish them from $x, y, z$ in
$B^{\pi_2}$. Thus $\pi_4^2=\pi_2$, $B^{\pi_4}=(B^{\pi_2})^\theta$,
where $\theta(x)=-x$, $\theta(y)=-y$ and $\theta(z)=z$, and
$B^{\pi_4}$ is generated by $u=x^2/8=x_1^4/16$, $v=y^2/8=x_2^4/16$
and $w=z/2=x_1x_2/4$, subject to the relation $uv=w^4$.

Let $d\geq 1$. By \ref{fdsPmod} and the representation theory of
$sl_2$, the unique $d$-dimensional simple Poisson $B^{\pi_2}$-module
$M$ has basis $v_j$, $0\leq j\leq d-1$, such that
 \begin{gather*} \{x,v_j\}_M=j(d-j)v_{j-1},\quad \{z,v_j\}_M=(2j+1-d)v_j/2,\\
\{y,v_j\}_M=v_{j+1}\mbox{ if }j<d-1,\quad \{y,v_{d-1}\}_M=0.
\end{gather*}
As we have seen above, all finite-dimensional simple Poisson modules
over $B^{\pi_4}$ are one-dimensional so, unless $d=1$, $M$ cannot
remain simple over $B^{\pi_4}$. If
$J=xB^{\pi_2}+yB^{\pi_2}+zB^{\pi_2}$, the unique Poisson maximal
ideal of $B^{\pi_2}$, then $u,v\in J^2$ so, as $J^2\subseteq
\Pann_A(M)$, $\{u,v_j\}_M=\{v,v_j\}_M=0$ while
$\{w,v_j\}_M=(2j+1-d)v_j/4$ so $M$ splits into the direct sum of $d$
one-dimensional simple Poisson $B^{\pi_4}$-modules. This contrasts
with the situation for associative $\C$-algebras as covered by
\cite[Corollary 26.13(iv)]{passman}.
\end{example}
\begin{rmk}
The finite-dimensional simple Poisson modules for the other types of
Kleinian singularities are all one-dimensional. If
$f=x^2+y^2z+z^{n-1}$, $n\geq 4$ (type $D_n$), or $f=x^2+y^3+z^4$
(type $E_6$), or $f=x^2+y^2+yz^3$ (type $E_7$) or $f=x^2+y^3+z^5$
(type $E_8$), then, as for type $A_n$ with $n>1$, the Lie algebra
$\g(J)$ is solvable for the unique Poisson maximal ideal
$J=xA+yA+zA$.
\end{rmk}

\subsection{Further invariants of the torus}
\label{fourfour} We close Section~4 by observing some very different
behaviour, to that observed in Example~\ref{KleinianA3up}, in the
finite-dimensional simple Poisson modules in the passage from a
Poisson algebra $C$ to the Poisson algebra $C^G$ of invariants for
the action of a finite group $G$ of Poisson automorphisms. Let
$B=\C[x_1^{\pm 1},x_2^{\pm 1}]$, let $\pi$ be as in
Subsection~\ref{torus} and let $C=B^\pi$. Thus $C$ is generated by
$x=x_1+x_1^{-1}, y=x_2+x_2^{-1}$ and $z=x_1x_2+x_1^{-1}x_2^{-1}$,
subject to the relation $f=0$, where $f=xyz-x^2-y^2-z^2+4$, and the
Poisson bracket on $C$ is induced by the exact bracket $\{-,-\}_f$
on $A$. Let $\theta$ denote both the Poisson automorphism of $B$
such that $\theta(x_1)=x_1$ and $\theta(x_2)=-x_2$ and its
restriction to $C$. Thus $\theta(x)=x$, $\theta(y)=-y$ and
$\theta(z)=-z$. If $K$ is the Klein $4$-group generated by $\pi$ and
$\theta$ then $B^K=C^\theta$. It is a routine matter to check that
$C^\theta$ is generated by $x$, $v:=\frac{1}{2}y^2$,
$w:=\frac{1}{2}yz$ and $t:=\frac{1}{2}z^2$ subject to the relations
$w^2=vt$ and $2xw-x^2-2v-2t+4=0$, the latter being the expression of
the relation $f=0$ in terms of $x,v,w$ and $t$. Hence $C^{\theta}$
is generated by $x$, $v$ and $w$ subject to the relation $g=0$,
where $g=2xvw-x^2v-2v^2-2w^2+4v$. In $C^\theta$, $\{x,v\}=2xv-4w$,
$\{v,w\}=2vw-2vx$ and $\{w,x\}=2t-2v=2xw-x^2-4v+4$ so the Poisson
bracket is induced by the exact Poisson bracket $\{-,-\}_g$ on
$A^\prime:=\C[x,v,w]$. There are four Poisson maximal ideals of
$A^\prime$ containing $g$:
\begin{gather*}
I_1=(x-2)A^\prime+(v-2)A^\prime+(w-2)A^\prime,\; L_1=(x-2)A^\prime+vA^\prime+wA^\prime,\\
I_2=(x+2)A^\prime+(v-2)A^\prime+(w+2)A^\prime,\;
L_2=(x+2)A^\prime+vA^\prime+wA^\prime.
\end{gather*}
If $J$ is any one of these, then $g^2\in J$ and $\g(J)\simeq sl_2$.
Consequently, as was the case for $C=B^\pi$ in
Subsection~\ref{torus}, $C^\theta$ is $4$-homogeneous. However the
correspondence between finite-dimensional simple Poisson modules for
$C$ and those for $C^\theta$ is more subtle than one might expect.
With the $J_i$'s as in Subsection~\ref{torus}, $J_2\cap
C^\theta=J_3\cap C^\theta=I_1$ and $J_4\cap C^\theta=J_5\cap
C^\theta=I_2$. Let $M$ be a finite-dimensional simple Poisson
$C$-module annihilated by $J_2$, let $0\neq N$ be a Poisson
$C^\theta$-submodule of $M$ and let $n\in N$. As $M$ is annihilated
by the maximal ideal $J_2$ of $C$ every $\C$-subspace of $M$ is a
$C$-submodule. In particular $(y-2)N=0=(z-2)N$. Hence
\begin{gather*}
\{y,n\}_M=\tfrac{1}{2}y\{y,n\}_M=\tfrac{1}{2}\{v,n\}_M\in N \mbox{
and}\\
\{y,n\}_M+\{z,n\}_M=\tfrac{1}{2}z\{y,n\}_M+\tfrac{1}{2}y\{z,n\}_M=\{w,n\}_M\in
N.\end{gather*} Therefore $N$ is a Poisson $C$-submodule of $M$,
$N=M$ and $M$ is simple as a Poisson $C^\theta$-module. As $J_2\cap
C^\theta=I_1$, $M$ is isomorphic, as a Poisson $C^\theta$-module, to
the unique $(\dim M)$-dimensional simple Poisson $C^\theta$-module.
By the same argument, the same is true for any finite-dimensional
simple Poisson $C$-module annihilated by $J_3$. Similarly, if
$M^\prime$ and $M^{\prime\prime}$ are simple Poisson $C$-modules, of
the same finite dimension, annihilated by $J_4$ and $J_5$
respectively then, as Poisson $C^\theta$-modules, they are simple,
 isomorphic and annihilated by $I_2$.

On the other hand, the finite-dimensional simple Poisson
$C^\theta$-modules annihilated by $L_i$, $i=1,2$, are not Poisson
$C^\theta$-submodules of any simple Poisson $C$-modules. The unique
maximal ideal of $C$ containing $L_i$ is $(x-(6-4i))C+yC+zC$ which
is not Poisson.

There is a Poisson subalgebra $D$, of $C^\theta$, that is isomorphic
to $C$ and is such that the extensions $D\subset C^\theta$ and
$C^\theta\subset C$ exhibit the same behaviour of finite-dimensional
simple Poisson modules. If $\phi$ is the Poisson automorphism of $B$
with $\phi(x_1)=-x_1$ and $\phi(x_2)=x_2$ then $\phi$ restricts to
Poisson automorphisms of $C$, with $\phi(x)=-x$, $\phi(y)=y$ and
$\phi(z)=-z$, and $C^\theta$, with $\phi(x)=-x, \phi(v)=v$ and
$\phi(w)=-w$. Let $D=C^\phi$. It can be checked that $D$ is
generated by $u:=\frac{1}{2}x^2$, $v$ and $a:=\frac{1}{2}xw$ subject
to the relation $h=0$ where $h=2auv-2a^2-u^2v-uv^2+2uv$ and that its
Poisson bracket is induced by the exact bracket $\{-,-\}_{2h}$ on
$\C[a,u,v]$. There are four Poisson maximals:
\begin{gather*}
I_3=uD+aD+vD,\quad L_3=(u-2)D+aD+vD,\\
I_4=uD+aD+(v-2)D,\quad L_4=(u-2)D+(a-2)D+(v-2)D. \end{gather*} For
$d\geq 1$, the $d$-dimensional simple $\C^\theta$-modules
annihilated by $L_1$ or $L_2$, respectively $I_1$ or $I_2$, become
isomorphic simple Poisson $D$-modules annihilated by $L_4$,
respectively $L_3$, and there are two $d$-dimensional simple Poisson
$D$-modules, annihilated by $I_3$ or $I_4$, that are not
restrictions of simple Poisson $C^\theta$-modules. To see that
$D\simeq C$, let $G=\langle \pi, \theta, \phi \rangle$, which is
abelian, of order $8$, and let $H=\langle \theta, \phi \rangle$.
Then $B^H=\C[x_1^{\pm 2},x_2^{\pm 2}]\simeq B$ and
$D=(C^\theta)^\phi=((B^\pi)^\theta)^\phi=B^G=(B^H)^\pi\simeq
B^\pi=C$.

\section{The quotient varieties $\mathfrak{h}\oplus\mathfrak{h}^*/W.$}
Here we consider three Poisson algebras, computed by Alev and Foissy
\cite{alevfoissy}, that arise in the following way. Let $\g$ be a
simple Lie algebra, $\mathfrak{h}$ be a Cartan subalgebra and $W$ be
the Weyl group of $\g$. There is a diagonal action of $W$ on the
symplectic space $V:=\mathfrak{h}\oplus\mathfrak{h}^*$, and on the
symmetric algebra $S=S(V)$.  Alev and Foissy have computed the
Poisson homology in degree $0$ for the Poisson algebra $S^W$ for the
three types $A_2, B_2$ and $G_2$, where $\dim \mathfrak{h}=2$.

\begin{example}
In the $A_2$ case, $W$ is the symmetric group $S_3$ which acts as
permutations  on the Poisson algebra
$S^\prime=\C[x_1,x_2,x_3,y_1,y_2,y_3]$, with
$\{x_i,y_j\}=\delta_{ij}, \{x_i,x_j\}=0=\{y_i,y_j\}$, and on the
underlying $6$-dimensional symplectic space. The algebra
$S=S(\mathfrak{h}\oplus\mathfrak{h}^*)$ is generated by $x_1-x_2$,
$x_1-x_3$, $y_1-y_2$ and $y_1-y_3$. Set $a_1=2x_1-x_2-x_3$,
$a_2=-x_1+2x_2-x_3$, $a_3=-x_1-x_2+2x_3=-a_1-a_2$,
$b_1=2y_1-y_2-y_3$, $b_2=-y_1+2y_2-y_3$,
 and $b_3=-y_1-y_2+2y_3=-b_1-b_2$. Then
$W$ acts, as permutations, on $\{a_1,a_2,a_3\}$ and on
$\{b_1,b_2,b_3\}$ and $S$ is the polynomial algebra
$\C[a_1,a_2,b_1,b_2]$ with the Poisson bracket such that $
\{a_i,b_i\}=6$, $\{a_i,a_j\}=0=\{b_i,b_j\}$ and
$\{a_1,b_2\}=-3=\{a_2,b_1\}$. Alev and Foissy \cite[Th\'{e}or\`{e}me
13]{alevfoissy} identify the following $W$-invariants in $S$:
\begin{gather*}
g_1=(a_1^2+a_2^2+a_1a_3)/9,\quad
g_2=(b_1^2+b_2^2+b_1b_3)/9,\\
g_3=-(2a_1b_1+a_1b_2+a_2b_1+2a_2b_2)/9,\\
m_1=(a_1a_2^2+a_2a_1^2)/9,\quad
m_2=(b_1b_2^2+b_2b_1^2)/9,\\
m_3=(2a_1b_1b_2+2a_2b_1b_2+a_1b_2^2+a_2b_1^2)/9,\\
m_4=(2b_1a_1a_2+2b_2a_1a_2+b_1a_2^2+b_2a_1^2)/9,
\end{gather*}
and that these generate $S^W$ subject to the relations $r_i=0$,
where, for $1\leq i\leq 5$,
\begin{align*}
r_1=&g_3m_3+3g_1m_2+g_2m_4,\\ r_2=&g_3m_4+3g_2m_1+g_1m_3,\\
r_3=&m_3^2-12g_1g_2^2+3g_2g_3^2-3m_2m_4,\\
r_4=&m_4^2-12g_2g_1^2+3g_1g_3^1-3m_1m_3,\\
r_5=&m_3m_4-9m_1m_2+12g_1g_2g_3+3g_3^3.
\end{align*}
In \cite[Th\'{e}or\`{e}me 13]{alevfoissy}, there is a sixth relation
which can be deduced from the overlap ambiguity
$(g_3m_3)m_4=g_3(m_3m_4)$ between the first terms of $r_1$ and $r_5$
and which is used, with the other relations, to present $S^W$ as a
free module over the polynomial algebra $\C[g_1,m_1,g_2,m_2]$ with
basis $\{1,g_3,g_3^2,g_3^3,m_3,m_4\}$. In \cite{alevfoissy},
$g_1,g_2,g_3,m_1,m_2,m_3,m_4$ are written as
$S_1,S_2,H,T_1,T_2,U_1,U_2$ respectively. The reason for our change
will emerge shortly.

A table showing the Poisson bracket on $S^W$ in terms of these
generators is given in \cite[6.1.5]{alevfoissy} and it is observed
in \cite[6.1.6]{alevfoissy} that this lifts to the polynomial
algebra $P_7:=\C[g_1,g_2,g_3,m_1,m_2,m_3,m_4]$. In $P_7$, the
elements $r_i$, $1\leq i\leq 5$, generate a Poisson ideal $I$ such
that, as Poisson algebras, $S^W\simeq P_7/I$.

From \cite[6.1.5]{alevfoissy}, $g_1=\frac{1}{2}\{g_3,g_1\}$,
$g_2=\frac{1}{2}\{g_2,g_3\}$, $g_3=\{g_2,g_1\}$,
$m_1=\frac{1}{3}\{g_3,m_1\}$, $m_2=\frac{1}{3}\{m_2,g_3\}$,
$m_3=\{g_1,m_2\}$ and $m_4=\{m_1,g_2\}$. So if $V$ is the generating
subspace, spanned by the $g_i$'s and $m_j$'s, then $V\subseteq
\{V,V\}$. As $\{V,V\}\subseteq K$ for any Poisson maximal ideal $K$
of $P_7$, it follows that $J:=(g_1,g_2,g_3,m_1,m_2,m_3,m_4)$ is the
unique Poisson maximal ideal of $P_7$ and that $\J:=J/I$ is the
unique Poisson maximal ideal of $S^W$. As $r_i\in J^2$ for $1\leq
i\leq 5$, $\g(J)\simeq \g(\J)$ so $P_7$ and $S^W$ share their
finite-dimensional simple Poisson modules.

There is a Poisson $\Z$-grading of $S^\prime=\oplus S^\prime_i$, in
the sense that $S^\prime_iS^\prime_j\subseteq S^\prime_{i+j}$ and
$\{S^\prime_i,S^\prime_j\}\subseteq S^\prime_{i+j}$, in which, for
each $i$, $a_i\in S^\prime_1$ and $b_i\in S^\prime_{-1}$. This is
inherited by $S^W$ and can be lifted to $P_7$, where the generators
$g_1,g_2,g_3,m_1,m_2,m_3,m_4$ have degrees $2,-2,0,3,-3,-1,1$
respectively. The table in \cite[6.1.5]{alevfoissy} shows that
$\{g_3,p\}=\deg(p)p$ if $p$ is any of these generators and, by the
derivation properties of $\{-,-\}$, the same conclusion follows for
any $p\in P_7$. Deleting quadratic terms from the table in
\cite[6.1.5]{alevfoissy}, we find that, in $\g(J)$,
\begin{gather*}
[g_1,g_2]=-g_3,\;[g_2,g_3]=2g_2,\;[g_1,g_3]=-2g_1,\\
[m_i,m_j]=0, 1\leq i,j\leq 4,\\
[g_1,m_1]=0,\;[g_1,m_2]=m_3,\;[g_1,m_3]=2m_4,\;[g_1,m_4]=3m_1,\\
[g_2,m_1]=-m_4,\;[g_2,m_2]=0,\;[g_2,m_3]=-3m_2,\;[g_2,m_4]=-2m_3,\\
[g_3,m_1]=3m_1,\;[g_3,m_2]=-3m_2,\;[g_3,m_3]=-m_3,\;[g_3,m_4]=m_4.
\end{gather*}
From this it follows that
\begin{itemize}
\item $[\g(J),\g(J)]=\g(J)$;
\item $g_1,g_2$ and $g_3$ generate a subalgebra $\s$ isomorphic to $sl_2$;
\item $m_1,m_2,m_3$ and $m_4$ generate an abelian subalgebra $\m$;
\item $\m$ is an ideal of $\g$ and is the radical of
$\g(J)$;
\item by the action of $\ad g_3$, every non-zero Poisson
$\s$-submodule of $\m$ contains  $m_i$ for some $i$; \item $\m$ is a
simple $\s$-module.
\end{itemize}
Thus $\g(J)$ is isomorphic to the Lie algebra obtained by taking the
direct sum of $sl_2$ and its unique $4$-dimensional simple module
$M$ and extending the Lie bracket from $sl_2$ by taking
$[s,m]=[s,m]_M$ and $[m,n]=0$ for $s\in sl_2$ and $m,n\in M$. By
\cite[Proposition 1.7.1(i)]{dix}, every finite-dimensional simple
$\g(J)$-module is annihilated by $\m$ and becomes a simple
$sl_2$-module. By Theorem~\ref{fdsPmod}, $S^W$ and $P_7$ are both
 $1$-homogeneous.

The Poisson algebra $P_7$ is intriguing and it would be interesting
to identify a general class of Poisson brackets on polynomial
algebras of which this is an example. The following Proposition
gives some of its features. \label{Atwo}
\end{example}
\begin{prop}\label{idealsofP}
Let $J$ be the unique Poisson maximal ideal of $P_7$.

(i) For $1\leq i\leq 5$, $r_i$ is not Poisson central in $P_7$.

(ii) The ideals $I:=\sum_1^5 r_iP_7$ and $I_1:=\sum_3^5 r_iP_7+IJ$
are Poisson ideals of $P_7$. Let $M$ be the Poisson module $I/IJ$,
$N:=I_1/IJ$ and $N^\prime:=I/I_1$. Then $N$ is isomorphic to the
unique $3$-dimensional simple Poisson $P_7$-module and $N^\prime$ is
isomorphic to the unique $2$-dimensional simple Poisson
$P_7$-module. Moreover $N$ is the unique proper non-zero Poisson
submodule of $M$ and $M$ is not Poisson semisimple.
\end{prop}
\begin{proof}
It is a routine calculation to use the table in
\cite[6.1.5]{alevfoissy} to compute $\{g,r_i\}$, where $g$ is one of
the seven generators of $P_7$ and $1\leq i\leq 5$. We have seen that
$\{g_3,p\}=(\deg p) p$ for all $p\in P_7$. The elements
$r_1,r_2,r_3,r_4,r_5$ have degrees $-1,1,-2,2,0$ respectively so
$r_i$ is not Poisson central if $1\leq i\leq 4$. For the other six
generators, the non-zero brackets are displayed below. The
non-centrality of $r_5$ follows from its bracket with $g_1$.
\begin{gather*}
\{g_1,r_1\}=r_2,\; \{g_1,r_3\}=r_5,\;\{g_1,r_5\}=2r_4,\\
\{g_2,r_2\}=-r_1,\; \{g_2,r_4\}=-r_5,\;\{g_2,r_5\}=-2r_3,\\
\{m_1,r_1\}=r_4,\; \{m_1,r_3\}=-6(g_1r_1+g_3r_2),\;\{m_1,r_5\}=6g_1r_2,\\
\{m_2,r_2\}=-r_3,\; \{m_2,r_4\}=6(g_2r_2+g_3r_1),\;\{m_2,r_5\}=-6g_2r_1,\\
\{m_3,r_1\}=r_3,\; \{m_3,r_2\}=-r_5,\; \{m_3,r_3\}=6g_2r_1,\\
 \{m_3,r_4\}=-12g_1r_1,\;\{m_3,r_5\}=-6(g_3r_1-3g_2r_2),\\
\{m_4,r_1\}=r_5,\;\{m_4,r_2\}=-r_4,\; \{m_4,r_3\}=12g_2r_2,\\
\{m_4,r_4\}=-6g_1r_2,\; \{m_4,r_5\}=6(g_3r_2+3g_1r_1).
\end{gather*}
It is now apparent that $I$ and $I_1$ are Poisson ideals and that
$N$ is a Poisson submodule of $M$. The modules $M, N$ and $N^\prime$
have bases consisting of (the images of) $r_i$, where, respectively,
$1\leq i\leq 5$, $3\leq i\leq 5$ and $1\leq i\leq 2$.
 Every
finite-dimensional Poisson module $X$, over $P_7$ or $S^W$, is the
direct sum, over $\Z$, of the eigenspaces $E_\rho:=\{m\in
X:\{g_3,m\}_X=\rho m\}$. For $M$ these are $E_{-2}=\C r_3$,
$E_{-1}=\C r_1$, $E_{0}=\C r_5$, $E_{1}=\C r_2$ and $E_{2}=\C r_4$.
From the data displayed, the only sums of these that are Poisson
submodules are $0$, $M$ and $N$. Thus $N$ is simple and is the
unique Poisson submodule of $M$, which cannot then be Poisson
semisimple. The same method establishes the Poisson simplicity of
$N^\prime$, in which $\{g_1,r_1\}_{N^\prime}=r_2$,
$\{g_1,r_2\}_{N^\prime}=0=\{g_2,r_1\}_{N^\prime}$,
$\{g_2,r_2\}_{N^\prime}=-r_1$, $\{g_3,r_1\}_{N^\prime}=-r_1$ and
$\{g_3,r_2\}_{N^\prime}=r_2$.
\end{proof}

\begin{rmk}
Let $\s=\C g_1+\C g_2+\C g_3$, a Lie subalgebra of $P_7$ isomorphic
to $sl_2$. The Poisson $P_7$-module $M=I/IJ$ and its submodule
$N=I_1/IJ$ are $\s$-modules. In accordance with the representation
theory of $sl_2$, there is a simple $\s$-submodule $N^\prime$,
spanned by (the images of) $r_1$ and $r_2$, such that $M=N\oplus
N^\prime$. As we have seen in Proposition~\ref{idealsofP}, $M$ is
not semisimple as a Poisson $P_7$-module or as a Poisson
$S^W$-module. As $S$ itself is Poisson simple, the property that
every finite-dimensional Poisson module is semisimple is not
preserved by the passage from a Poisson algebra to the Poisson
algebra of invariants for the action of a finite group of Poisson
automorphisms. This is in contrast to the situation for associative
$\C$-algebras, see \cite{KS} or \cite[p. 470]{jw}. We take the
opportunity to correct an oversight in the latter reference, where
the argument given should have been attributed to M. Lorenz.
\end{rmk}

\begin{example}
\label{Gtwo} Next we consider, in rather less detail, the case where
$\g$ is of type $G_2$, appealing to \cite{alevfoissy} for
identification of the invariants. Here $S$ is as in
Example~\ref{Atwo} and the Weyl group, which we denote $W^\prime$,
is the direct product of the Weyl group $W$ for $A_2$ and a cyclic
group of order $2$, whose generator $\pi$ is such that
$\pi(a_i)=-a_i$ and $\pi(b_i)=-b_i$ for $i=1,2$. Consequently
$S^{W^\prime}=(S^W)^\theta$. By \cite[Th\'{e}or\`{e}me
15]{alevfoissy}, $S^{W^\prime}$ is generated by $g_1,g_2,g_3,
n_1:=m_1^2, n_2:=m_2^2, n_3:=m_1m_2, n_4:=m_1m_3, n_5:=m_1m_4,
n_6=m_2m_3$ and $n_7=m_2m_4$ subject to sixteen relations. The
Poisson bracket again lifts to the appropriate polynomial ring,
$P_{10}=\C[g_i,n_j,1\leq i\leq 3, 1\leq j\leq 7]$ which has a unique
Poisson maximal ideal $J=\sum_{i=1}^3 g_iP_{10}+\sum_{j=1}^7
n_jP_{10}$. As in Example~\ref{Atwo}, each of the relations has the
form $r_i=0$, where $r_i\in J^2$, so $P_{10}$ and $S^{W^\prime}$
share their finite-dimensional simple Poisson modules. The brackets
$\{g,h\}$ between generators can be computed from the data in
Example~\ref{Atwo} and it can then be seen, much as in that example,
that
 $[\g(J),\g(J)]=\g(J)$, that
$g_1,g_2$ and $g_3$ generate a subalgebra $\s$ isomorphic to $sl_2$,
that
 $n_1,n_2,n_3,n_4,n_5,n_6$ and $n_7$ generate an abelian subalgebra $\mathfrak{n}$, which  is the radical of $\g(J)$;
and is simple as an  $\s$-module, that $\g(J)$ is isomorphic to the
Lie algebra obtained by taking the direct sum of $sl_2$ and its
unique $7$-dimensional simple module $N$ and extending the Lie
bracket from $sl_2$ as in Example~\ref{Atwo} and that $S^W$ and
$P_{10}$ are both
 $1$-homogeneous Poisson algebras.
\end{example}
\begin{example}
Finally we consider the case where $\g$ is of type $B_2$. The
discussion is again based on the calculations of invariants in
\cite{alevfoissy}. Here
$S=S(\mathfrak{h}\oplus\mathfrak{h}^*)=\C[x_1,x_2,y_1,y_2]$,
$\{x_i,x_j\}=0=\{y_i,y_j\}$ and $\{x_i,y_j\}=\delta_{ij}$.  The Weyl
group $W$ is isomorphic to the dihedral group of order $8$ and is
generated by two elements that act on $S$ by $x_1\mapsto -x_1,
y_1\mapsto -y_1, x_2\mapsto x_2, y_2\mapsto y_2$ and by $x_1\mapsto
x_2, y_1\mapsto y_2, x_2\mapsto x_1, y_2\mapsto y_1$. The following
eight generators for $S^W$ are identified in \cite[Th\'{e}or\`{e}me
15]{alevfoissy}:
\begin{gather*}
g_1=x_1^2+x_2^2,\quad  g_2=y_1^2+y_2^2,\quad g_3=x_1y_1+x_2y_2,
\\
m_1=x_1^2x_2^2,\quad m_2=y_1^2y_2^2,\quad m_3=x_1y_1+x_2y_2,\\
m_4=x_1y_1^3+x_2y_2^3,\quad m_5=x_1^3y_1+x_2^3y_2.
\end{gather*}
As in Example~\ref{Atwo}, the Poisson bracket lifts from $S^W$ to a
Poisson bracket on the polynomial algebra $P_8=\C[g_i,m_j,1\leq
i\leq 3, 1\leq j\leq 5]$. It is a routine matter to compute
$\{g,h\}$ for each pair $g,h$ of generators and to conclude that
$P_8$ has a unique Poisson maximal ideal $J=\sum_{i=1}^3
g_iP_8+\sum_{j=1}^5 m_jP_8$. A set of ten defining relations,
$r_i=0$, for $S_W$ is given in \cite[Th\'{e}or\`{e}me
15]{alevfoissy} and each $r_i\in J^2$. Consequently $S^W$ has a
unique Poisson maximal ideal $\ov J$,  $\g(J)\simeq\g(\J)$  and
$P_8$ and $S^W$ share their finite-dimensional simple Poisson
modules. In $\g(J)$,
\begin{gather*}
[g_1,g_2]=4g_3,\; [g_1,g_3]=2g_1,\; [g_2,g_3]=-2g_2,\\
[m_j,m_k]=0\mbox{ for }1\leq j,k\leq 5,\\
[g_1,m_1]=[g_2,m_2]=[g_3,m_3]=0,\\ [g_1,m_2]=-4m_4, [g_1,m_3]=-2m_5,
[g_1,m_4]=-12m_3, [g_1,m_5]=-4m_1,\\
[g_2,m_1]=4m_5,  [g_2,m_3]=2m_4,
[g_2,m_4]=4m_2, [g_2,m_5]=12m_3,\\
[g_3,m_1]=-4m_1, [g_3,m_2]=4m_2,  [g_3,m_4]=2m_4, [g_3,m_5]=-2m_5.
\end{gather*}
From this it follows, much as in Example~\ref{Atwo}, that
$[\g(J),\g(J)]=\g(J)$, that $g_1,g_2$ and $g_3$ generate a
subalgebra $\s$ isomorphic to $sl_2$, that $m_1,m_2,m_3,m_4$ and
$m_5$ generate an abelian subalgebra $\m$, which is the radical of
$\g(J)$ and is simple as an $\s$-module, that $\g(J)$ is isomorphic
to the Lie algebra obtained by taking the direct sum of $sl_2$ and
its unique $5$-dimensional simple module $M$ and extending the Lie
bracket from $sl_2$ as before and that $S^W$ and $P_8$ are
 $1$-homogeneous Poisson algebras.
\label{Btwo}
\end{example}

\end{document}